\documentclass[12pt]{article}

\usepackage{amsmath, epsfig, cite}
\usepackage{amssymb}
\usepackage{amsfonts}
\usepackage{latexsym}
\usepackage{amsthm}

\newtheorem{thm}{Theorem}[section]

\numberwithin{equation}{section}

\renewcommand{\thefootnote}

\begin{document}

\begin{center}
{\large\bf On some conjectures of Z.-W. Sun involving harmonic
numbers
 \footnote{ The work is supported by the National Natural Science Foundation of China (No. 12071103).}}
\end{center}

\renewcommand{\thefootnote}{$\dagger$}

\vskip 2mm \centerline{Chuanan Wei}
\begin{center}
{School of Biomedical Information and Engineering\\ Hainan Medical
University, Haikou 571199, China
\\
 Email address: weichuanan78@163.com}
\end{center}


\vskip 0.7cm \noindent{\bf Abstract.} Harmonic numbers are
significant in various branches of number theory.  With the help of
the digamma function, we prove ten conjectural series of Z.-W. Sun
involving harmonic numbers. Several ones of them are also series
expansions of $\log2/\pi^2$.



\section{Introduction}

 For a complex number $x$, define the generalized harmonic numbers to
be
\[H_{0}(x)=0
\quad\text{and}\quad H_{n}(x)
  =\sum_{k=1}^n\frac{1}{x+k}.\]
When $x=0$, they reduce to classical harmonic numbers:
\[H_{0}=0\quad\text{and}\quad
H_{n}=\sum_{k=1}^n\frac{1}{k}.\] For a nonnegative integer $m$,
define the shifted-factorial by
\begin{align*}
(x)_0=1 \quad\text{and}\quad (x)_m=x(x+1)\cdots(x+m-1) \quad
\text{when} \quad m\in \mathbb{Z}^{+}.
\end{align*}
 For a differentiable function $f(x)$, define the derivative operator
$\mathcal{D}_x$ as
\begin{align*}
\mathcal{D}_xf(x)=\frac{d}{dx}f(x).
 \end{align*}
 Then it is routine to show that
$$\mathcal{D}_x\:(1+x)_r=(1+x)_rH_r(x),$$
where $r$ is an arbitrary nonnegative integer. Several nice harmonic
number identities from differentiation of the shifted-factorials can
be seen in the papers \cite{Paule,Sofo,Wang-Wei}.

There exist a lot of interesting $\pi$-formulas in the literature.
Four series for $1/\pi^2$ due to Guillera
\cite{Guillera-a,Guillera-b,Guillera-c} can be laid out as follows:
\begin{align}
&\sum_{k=0}^{\infty}(820k^2+180k+13)\frac{\binom{2k}{k}^5}{(-2^{20})^k}=\frac{128}{\pi^2},
\label{Guillera-a}\\[1mm]
&\quad\sum_{k=0}^{\infty}(20k^2+8k+1)\frac{\binom{2k}{k}^5}{(-2^{12})^k}=\frac{8}{\pi^2},
\label{Guillera-b}\\[1mm]
&\:\:\sum_{k=0}^{\infty}(120k^2+34k+3)\frac{\binom{2k}{k}^4\binom{4k}{2k}}{2^{16k}}=\frac{32}{\pi^2}.
\label{Guillera-c}\\[1mm]
&\:\:\:\sum_{k=0}^{\infty}(74k^2+27k+3)\frac{\binom{2k}{k}^4\binom{3k}{k}}{2^{12k}}=\frac{48}{\pi^2}.
\label{Guillera-d}
\end{align}
For more conclusions on $\pi$-formulas, the reader is referred to
 the papers \cite{Au,Guo,Wang}.

On the basis of the three identities from \emph{Mathematica}:
\begin{align*}
&\:\:\sum_{k=0}^{\infty}\frac{\binom{2k}{k}\binom{3k}{k}}{54^k}=\frac{3\Gamma(\frac{1}{3})^3}{4\pi^2\sqrt[3]{2}},
\\[1mm]
&\sum_{k=0}^{\infty}\frac{\binom{2k}{k}\binom{4k}{2k}}{128^k}=\frac{\sqrt{\pi}}{\Gamma(\frac{5}{8})\Gamma(\frac{7}{8})},
\\[1mm]
&\sum_{k=0}^{\infty}\frac{\binom{3k}{k}\binom{6k}{3k}}{864^k}=\frac{\sqrt{\pi}}{\Gamma(\frac{7}{12})\Gamma(\frac{11}{12})},
\end{align*}
Sun \cite[Equations (2.9), (2.11) and (2.16)]{Sun} proposed the
following three conjectures.

\begin{thm}\label{theorem-a}
\begin{align}
\sum_{k=0}^{\infty}\frac{\binom{2k}{k}\binom{3k}{k}}{54^k}(3H_{3k}-H_k)=\frac{(3\log2)\Gamma(\frac{1}{3})^3}{4\pi^2\sqrt[3]{2}}.
 \label{wei-a}
\end{align}
\end{thm}

\begin{thm}\label{theorem-b}
\begin{align}
\sum_{k=0}^{\infty}\frac{\binom{2k}{k}\binom{4k}{2k}}{128^k}(2H_{4k}-H_{2k})=\frac{(\log2)\sqrt{\pi}}{2\Gamma(\frac{5}{8})\Gamma(\frac{7}{8})}.
 \label{wei-b}
\end{align}
\end{thm}

\begin{thm}\label{theorem-c}
\begin{align}
\sum_{k=0}^{\infty}\frac{\binom{3k}{k}\binom{6k}{3k}}{864^k}(6H_{6k}-3H_{3k}-2H_{2k}+H_k)=\frac{(\log2)\sqrt{\pi}}{\Gamma(\frac{7}{12})\Gamma(\frac{11}{12})}.
 \label{wei-c}
\end{align}
\end{thm}

Motivated by \eqref{Guillera-a} and \eqref{Guillera-b}, Sun
\cite[Equations (4.22), (4.15) and (4.16)]{Sun} proposed the
following three conjectures.

\begin{thm}\label{thm-a}
\begin{align}
\sum_{k=0}^{\infty}\frac{\binom{2k}{k}^5}{(-2^{20})^k}\Big\{(820k^2+180k+13)(H_{2k}-H_k)+164k+18\Big\}=\frac{256\log2}{\pi^2}.
 \label{eq:wei-a}
\end{align}
\end{thm}

\begin{thm}\label{thm-b}
\begin{align}
\sum_{k=0}^{\infty}\frac{\binom{2k}{k}^5}{(-2^{12})^k}\Big\{(20k^2+8k+1)H_k-6k-1\Big\}=-\frac{16\log2}{\pi^2}.
 \label{eq:wei-b}
\end{align}
\end{thm}

\begin{thm}\label{thm-c}
\begin{align}
\sum_{k=0}^{\infty}\frac{\binom{2k}{k}^5}{(-2^{12})^k}\Big\{5(20k^2+8k+1)H_{2k}-10k-1\Big\}=-\frac{32\log2}{\pi^2}.
 \label{eq:wei-c}
\end{align}
\end{thm}

Encouraged by \eqref{Guillera-c} and \eqref{Guillera-d},  Sun
\cite[Equations (4.34), (4.33), (4.27) and (4.28)]{Sun} suggested
the following four conjectures.

\begin{thm}\label{thm-d}
\begin{align}
\sum_{k=0}^{\infty}\frac{\binom{2k}{k}^4\binom{4k}{2k}}{2^{16k}}\Big\{(120k^2+34k+3)(H_{2k}-2H_k)+68k+9\Big\}=\frac{128\log2}{\pi^2}.
 \label{eq:wei-d}
\end{align}
\end{thm}

\begin{thm}\label{thm-e}
\begin{align}
\sum_{k=0}^{\infty}\frac{\binom{2k}{k}^4\binom{4k}{2k}}{2^{16k}}\Big\{2(120k^2+34k+3)H_{4k}-16k-1\Big\}=0.
 \label{eq:wei-e}
\end{align}
\end{thm}

\begin{thm}\label{thm-f}
\begin{align}
\sum_{k=0}^{\infty}\frac{\binom{2k}{k}^4\binom{3k}{k}}{2^{12k}}\Big\{(74k^2+27k+3)H_{2k}-17k-3\Big\}=0.
 \label{eq:wei-f}
\end{align}
\end{thm}

\begin{thm}\label{thm-g}
\begin{align}
\sum_{k=0}^{\infty}\frac{\binom{2k}{k}^4\binom{3k}{k}}{2^{12k}}\Big\{(74k^2+27k+3)(51H_{3k}+250H_{2k}-153H_k)+15\Big\}=\frac{9792\log2}{\pi^2}.
 \label{eq:wei-g}
\end{align}
\end{thm}

Define the digamma function $\psi(x)$ to be
\begin{align*}
\psi(x)=\frac{d}{dx}\big\{\log\Gamma(x)\big\},
\end{align*}
where $\Gamma(x)$ is the familiar gamma function. Some special
values of the digamma function (cf. \cite{Liu}) are known. For
example,
\begin{align}
&\qquad\quad\psi(1)=-\gamma,
 \label{digamma-a}\\[1mm]
&\quad\psi(\tfrac{1}{2})=-\gamma-2\log2,
 \label{digamma-b}
\end{align}
where $\gamma$ is the Euler-Mascheroni constant.

The structure of the paper is organized as follows. According to the
properties of the digamma function, we shall provide the proof of
Theorems \ref{theorem-a}-\ref{theorem-c} in Section 2. Similarly,
the proof of Theorems \ref{thm-a}-\ref{thm-c} and Theorems
\ref{thm-d}-\ref{thm-g} will respectively be displayed in Sections 3
and 4.

\section{Proof of Theorems \ref{theorem-a}-\ref{theorem-c}}
For the goal of proving Theorems \ref{theorem-a}-\ref{theorem-c}, we
require Bailey's $_2F_1(\frac{1}{2})$ summation formula (cf.
\cite[P.17]{Bailey}) and Gauss' $_2F_1(\frac{1}{2})$ summation
formula (cf. \cite[P.17]{Bailey}):
\begin{align}
&{_{2}F_{1}}\left[\begin{array}{cccccccc}
  a,1-a\\
 b
\end{array};\frac{1}{2}\right]
=\frac{\Gamma(\frac{b}{2})\Gamma(\frac{1+b}{2})}{\Gamma(\frac{a+b}{2})\Gamma(\frac{1-a+b}{2})},
 \label{eq:Bailey}
 \\[1mm]
&\:\:{_{2}F_{1}}\left[\begin{array}{cccccccc}
  a,b\\
 \frac{1+a+b}{2}
\end{array};\frac{1}{2}\right]
=\frac{\Gamma(\frac{1}{2})\Gamma(\frac{1+a+b}{2})}{\Gamma(\frac{1+a}{2})\Gamma(\frac{1+b}{2})},
 \label{eq:Gauss}
\end{align}
where the hypergeometric has been defined by
$$
_{r+1}F_{r}\left[\begin{array}{c}
a_1,a_2,\ldots,a_{r+1}\\
b_1,b_2,\ldots,b_{r}
\end{array};\, z
\right] =\sum_{k=0}^{\infty}\frac{(a_1)_k(a_2)_k\cdots(a_{r+1})_k}
{(1)_k(b_1)_k\cdots(b_{r})_k}z^k.
$$

Firstly, we shall prove Theorem \ref{theorem-a}.

\begin{proof}[{\bf{Proof of Theorem \ref{theorem-a}}}]
We comprehend that these series in \eqref{eq:Bailey} and
\eqref{eq:Gauss} are both uniformly convergent for $b\in
\mathbb{C}$. Apply $\mathcal{D}_{b}$ on both sides of
\eqref{eq:Bailey} and \eqref{eq:Gauss} to get
\begin{align}
&\sum_{k=0}^{\infty}\bigg(\frac{1}{2}\bigg)^{k-1}\frac{(a)_k(1-a)_k}{(1)_{k}(b)_{k}}H_k(b-1)=\frac{\Gamma(\frac{b}{2})\Gamma(\frac{1+b}{2})}{\Gamma(\frac{a+b}{2})\Gamma(\frac{1-a+b}{2})}
\notag\\[1mm]
&\:\:\times
\Big\{\psi(\tfrac{a+b}{2})+\psi(\tfrac{1-a+b}{2})-\psi(\tfrac{b}{2})-\psi(\tfrac{1+b}{2})\Big\},
 \label{eq:Bailey-a}
\\[1mm]
&\sum_{k=0}^{\infty}\bigg(\frac{1}{2}\bigg)^{k-1}\frac{(a)_k(b)_k}{(1)_{k}(\frac{1+a+b}{2})_{k}}\Big\{H_k(b-1)-\tfrac{1}{2}H_k(\tfrac{a+b-1}{2})\Big\}
\notag\\[1mm]
&\:\:=\frac{\Gamma(\frac{1}{2})\Gamma(\frac{1+a+b}{2})}{\Gamma(\frac{1+a}{2})\Gamma(\frac{1+b}{2})}
\Big\{\psi(\tfrac{1+a+b}{2})-\psi(\tfrac{1+b}{2})\Big\}.
 \label{eq:Gauss-a}
\end{align}
The $(a,b)=(\frac{1}{3},1)$ case of \eqref{eq:Bailey-a} reads
\begin{align}
\sum_{k=0}^{\infty}\frac{\binom{2k}{k}\binom{3k}{k}}{54^k}H_k=\frac{\Gamma(\frac{1}{2})}{2\Gamma(\frac{2}{3})\Gamma(\frac{5}{6})}
\Big\{\psi(\tfrac{2}{3})+\psi(\tfrac{5}{6})-\psi(\tfrac{1}{2})-\psi(1)\Big\}.
 \label{eq:Bailey-b}
\end{align}
The $(a,b)=(\frac{1}{3},\frac{2}{3})$  and
$(a,b)=(\frac{2}{3},\frac{1}{3})$ cases of \eqref{eq:Gauss-a} are
\begin{align}
\sum_{k=0}^{\infty}\frac{\binom{2k}{k}\binom{3k}{k}}{54^k}\Big\{H_k(-\tfrac{1}{3})-\tfrac{1}{2}H_k\Big\}
=\frac{\Gamma(\frac{1}{2})}{2\Gamma(\frac{2}{3})\Gamma(\frac{5}{6})}
\Big\{\psi(1)-\psi(\tfrac{5}{6})\Big\},
 \label{eq:Gauss-b}\\[1mm]
\sum_{k=0}^{\infty}\frac{\binom{2k}{k}\binom{3k}{k}}{54^k}\Big\{H_k(-\tfrac{2}{3})-\tfrac{1}{2}H_k\Big\}
=\frac{\Gamma(\frac{1}{2})}{2\Gamma(\frac{2}{3})\Gamma(\frac{5}{6})}
\Big\{\psi(1)-\psi(\tfrac{2}{3})\Big\}.
 \label{eq:Gauss-c}
\end{align}
By means of \eqref{digamma-a}, \eqref{digamma-b} and the relation
\begin{align*}
H_k(-\tfrac{1}{3})+H_k(-\tfrac{2}{3})=3H_{3k}-H_k,
\end{align*}
the sum of \eqref{eq:Bailey-b}, \eqref{eq:Gauss-b} and
\eqref{eq:Gauss-c} results in \eqref{wei-a}.
\end{proof}

Secondly, we prepare to prove Theorem \ref{theorem-b}.

\begin{proof}[{\bf{Proof of Theorem \ref{theorem-b}}}]
The $(a,b)=(\frac{1}{4},1)$ case of \eqref{eq:Bailey-a} yields
\begin{align}
\sum_{k=0}^{\infty}\frac{\binom{2k}{k}\binom{4k}{2k}}{128^k}H_k=\frac{\Gamma(\frac{1}{2})}{2\Gamma(\frac{5}{8})\Gamma(\frac{7}{8})}
\Big\{\psi(\tfrac{5}{8})+\psi(\tfrac{7}{8})-\psi(\tfrac{1}{2})-\psi(1)\Big\}.
 \label{eq:Bailey-c}
\end{align}
The $(a,b)=(\frac{1}{4},\frac{3}{4})$  and
$(a,b)=(\frac{3}{4},\frac{1}{4})$ cases of \eqref{eq:Gauss-a} are
\begin{align}
\sum_{k=0}^{\infty}\frac{\binom{2k}{k}\binom{4k}{2k}}{128^k}\Big\{H_k(-\tfrac{1}{4})-\tfrac{1}{2}H_k\Big\}
=\frac{\Gamma(\frac{1}{2})}{2\Gamma(\frac{5}{8})\Gamma(\frac{7}{8})}
\Big\{\psi(1)-\psi(\tfrac{7}{8})\Big\},
 \label{eq:Gauss-d}\\[1mm]
\sum_{k=0}^{\infty}\frac{\binom{2k}{k}\binom{4k}{2k}}{128^k}\Big\{H_k(-\tfrac{3}{4})-\tfrac{1}{2}H_k\Big\}
=\frac{\Gamma(\frac{1}{2})}{2\Gamma(\frac{5}{8})\Gamma(\frac{7}{8})}
\Big\{\psi(1)-\psi(\tfrac{5}{8})\Big\}.
 \label{eq:Gauss-e}
\end{align}
Via \eqref{digamma-a}, \eqref{digamma-b} and the relation
\begin{align*}
H_k(-\tfrac{1}{4})+H_k(-\tfrac{3}{4})=4H_{4k}-2H_{2k},
\end{align*}
the sum of \eqref{eq:Bailey-c}, \eqref{eq:Gauss-d} and
\eqref{eq:Gauss-e} leads to \eqref{wei-b}.
\end{proof}

Thirdly, we intend to prove Theorem \ref{theorem-c}.

\begin{proof}[{\bf{Proof of Theorem \ref{theorem-c}}}]
The $(a,b)=(\frac{1}{6},1)$ case of \eqref{eq:Bailey-a} reads
\begin{align}
\sum_{k=0}^{\infty}\frac{\binom{3k}{k}\binom{6k}{3k}}{864^k}H_k=\frac{\Gamma(\frac{1}{2})}{2\Gamma(\frac{7}{12})\Gamma(\frac{11}{12})}
\Big\{\psi(\tfrac{7}{12})+\psi(\tfrac{11}{12})-\psi(\tfrac{1}{2})-\psi(1)\Big\}.
 \label{eq:Bailey-d}
\end{align}
The $(a,b)=(\frac{1}{6},\frac{5}{6})$  and
$(a,b)=(\frac{5}{6},\frac{1}{6})$ cases of \eqref{eq:Gauss-a} yield
\begin{align}
\sum_{k=0}^{\infty}\frac{\binom{3k}{k}\binom{6k}{3k}}{864^k}\Big\{H_k(-\tfrac{1}{6})-\tfrac{1}{2}H_k\Big\}
=\frac{\Gamma(\frac{1}{2})}{2\Gamma(\frac{7}{12})\Gamma(\frac{11}{12})}
\Big\{\psi(1)-\psi(\tfrac{11}{12})\Big\},
 \label{eq:Gauss-f}\\[1mm]
\sum_{k=0}^{\infty}\frac{\binom{3k}{k}\binom{6k}{3k}}{864^k}\Big\{H_k(-\tfrac{5}{6})-\tfrac{1}{2}H_k\Big\}
=\frac{\Gamma(\frac{1}{2})}{2\Gamma(\frac{7}{12})\Gamma(\frac{11}{12})}
\Big\{\psi(1)-\psi(\tfrac{7}{12})\Big\}.
 \label{eq:Gauss-g}
\end{align}
Through \eqref{digamma-a}, \eqref{digamma-b} and the relation
\begin{align*}
H_k(-\tfrac{1}{6})+H_k(-\tfrac{5}{6})=6H_{6k}-3H_{3k}-2H_{2k}+H_k,
\end{align*}
the sum of \eqref{eq:Bailey-d}, \eqref{eq:Gauss-f} and
\eqref{eq:Gauss-g} gives \eqref{wei-c}.
\end{proof}
\section{Proof of Theorems \ref{thm-a}-\ref{thm-c}}
Firstly, we shall prove Theorem \ref{thm-a}.

\begin{proof}[{\bf{Proof of Theorem \ref{thm-a}}}]
A transformation formula for hypergeometric series (cf.
\cite[Theorem 31]{Chu-b}) may be written as
\begin{align}
&\sum_{k=0}^{\infty}(-1)^k\frac{(b)_k(c)_k(d)_k(e)_k(1+a-b-c)_k(1+a-b-d)_{k}(1+a-b-e)_{k}}{(1+a-b)_{2k}(1+a-c)_{2k}(1+a-d)_{2k}(1+a-e)_{2k}}
\notag\\[1mm]
&\quad\times\frac{(1+a-c-d)_k(1+a-c-e)_{k}(1+a-d-e)_{k}}{(1+2a-b-c-d-e)_{2k}}\alpha_k(a,b,c,d,e)
\notag\\[1mm]
&\:=\sum_{k=0}^{\infty}(a+2k)\frac{(b)_k(c)_k(d)_k(e)_k}{(1+a-b)_{k}(1+a-c)_{k}(1+a-d)_{k}(1+a-e)_{k}},
\label{equation-aa}
\end{align}
where $\mathfrak{R}(1+2a-b-c-d-e)>0$ and
\begin{align*}
&\alpha_k(a,b,c,d,e)\\[1mm]
&\:=\frac{(1+2a-b-c-d+3k)(a-e+2k)}{1+2a-b-c-d-e+2k}+\frac{(e+k)(1+a-b-c+k)}{(1+a-b+2k)(1+a-d+2k)}
\\[1mm]
&\quad\times\frac{(1+a-b-d+k)(1+a-c-d+k)(2+2a-b-d-e+3k)}{(1+2a-b-c-d-e+2k)(2+2a-b-c-d-e+2k)}
\\[1mm]
&\:+\frac{(c+k)(e+k)(1+a-b-c+k)(1+a-b-d+k)}{(1+a-b+2k)(1+a-c+2k)(1+a-d+2k)(1+a-e+2k)}
\\[1mm]
&\quad\times\frac{(1+a-b-e+k)(1+a-c-d+k)(1+a-d-e+k)}{(1+2a-b-c-d-e+2k)(2+2a-b-c-d-e+2k)}.
\end{align*}

Select $(a,b,c,e)=(\frac{1}{2},\frac{1}{2},\frac{1}{2},\frac{1}{2})$
in \eqref{equation-aa} to obtain
\begin{align}
\sum_{k=0}^{\infty}\bigg(\frac{-1}{64}\bigg)^k\frac{(\frac{1}{2})_k^3(d)_k(1-d)_k^3}{(1)_{k}^3(\frac{3}{2}-d)_{2k}^2}A_k(d)
=\frac{2}{1-2d}\:{_{5}F_{4}}\left[\begin{array}{cccccccc}
  \frac{1}{2},\frac{5}{4},\frac{1}{2},\frac{1}{2},d\\
  \frac{1}{4},1,1,\frac{3}{2}-d
\end{array};1\right],
\label{eq:wei-aa}
\end{align}
where
\begin{align*}
A_k(d) &=4k(1-d+3k)+\frac{2(1+2k)(2-d+3k)(1-d+k)^2}{(3-2d+4k)^2}
\\[1mm]
&\quad+\frac{(1+2k)(1-d+k)^3}{2(3-2d+4k)^2}.
\end{align*}
 Evaluating the $_5F_4$ series in
\eqref{eq:wei-aa} by Dougall's theorem (cf. \cite[P. 27]{Bailey}):
\begin{align}
&{_{5}F_{4}}\left[\begin{array}{cccccccc}
  a,1+\frac{a}{2},b,c,d\\
  \frac{a}{2},1+a-b,1+a-c,1+a-d
\end{array};1\right]
\notag\\[1mm]
&\:\:=
\frac{\Gamma(1+a-b)\Gamma(1+a-c)\Gamma(1+a-d)\Gamma(1+a-b-c-d)}{\Gamma(1+a)\Gamma(1+a-b-c)\Gamma(1+a-b-d)\Gamma(1+a-c-d)},
 \label{eq:Dougall}
\end{align}
 we have
\begin{align}
\sum_{k=0}^{\infty}\bigg(\frac{-1}{64}\bigg)^k\frac{(\frac{1}{2})_k^3(d)_k(1-d)_k^3}{(1)_{k}^3(\frac{3}{2}-d)_{2k}^2}A_k(d)
=\frac{2}{\pi}\frac{\Gamma(\frac{3}{2}-d)^2}{\Gamma(1-d)^2}.
\label{eq:wei-bb}
\end{align}

Notice that the series in \eqref{eq:wei-bb} is uniformly convergent
for $d\in \mathbb{C}$. Apply $\mathcal{D}_{d}$ on both sides of
\eqref{eq:wei-bb} to find
\begin{align*}
&\sum_{k=0}^{\infty}\bigg(\frac{-1}{64}\bigg)^k\frac{(\frac{1}{2})_k^3(d)_k(1-d)_k^3}{(1)_{k}^3(\frac{3}{2}-d)_{2k}^2}A_k(d)
\notag\\[1mm]
&\quad\times\Big\{H_{k}(d-1)-3H_{k}(-d)+2H_{2k}(\tfrac{1}{2}-d)\Big\}
\notag\\[1mm]
&\:+\sum_{k=0}^{\infty}\bigg(\frac{-1}{64}\bigg)^k\frac{(\frac{1}{2})_k^3(d)_k(1-d)_k^3}{(1)_{k}^3(\frac{3}{2}-d)_{2k}^2}\mathcal{D}_{d}A_k(d)
\notag\\[1mm]
&\:\:=\frac{4}{\pi}\frac{\Gamma(\frac{3}{2}-d)^2}{\Gamma(1-d)^2}
\Big\{\psi(1-d)-\psi(\tfrac{3}{2}-d)\Big\}.
\end{align*}
Setting $d=\frac{1}{2}$ in the last equation and using
\eqref{digamma-a} and \eqref{digamma-b}, we arrive at
\eqref{eq:wei-a} after some simplification.
\end{proof}

Secondly, we begin to prove Theorem \ref{thm-b}.

\begin{proof}[{\bf{Proof of Theorem \ref{thm-b}}}]
Another transformation formula for hypergeometric series (cf.
\cite[Theorem 9]{Chu-b}) may be expressed as
\begin{align}
&\sum_{k=0}^{\infty}\frac{(c)_k(d)_k(e)_k(1+a-b-c)_k(1+a-b-d)_{k}(1+a-b-e)_{k}}{(1+a-c)_{k}(1+a-d)_{k}(1+a-e)_{k}(1+2a-b-c-d-e)_{k}}
\notag\\[1mm]
&\quad\times\frac{(-1)^k}{(1+a-b)_{2k}}\beta_k(a,b,c,d,e)
\notag\\[1mm]
&\:=\sum_{k=0}^{\infty}(a+2k)\frac{(b)_k(c)_k(d)_k(e)_k}{(1+a-b)_{k}(1+a-c)_{k}(1+a-d)_{k}(1+a-e)_{k}},
\label{equation-bb}
\end{align}
where  $\mathfrak{R}(1+2a-b-c-d-e)>0$ and
\begin{align*}
\beta_k(a,b,c,d,e)&=\frac{(1+2a-b-c-d+2k)(a-e+k)}{1+2a-b-c-d-e+k}
\\[1mm]
&\quad+\frac{(1+a-b-c+k)(1+a-b-d+k)(e+k)}{(1+a-b+2k)(1+2a-b-c-d-e+k)}.
\end{align*}
Choose $(a,b,c,e)=(\frac{1}{2},\frac{1}{2},\frac{1}{2},\frac{1}{2})$
in \eqref{equation-bb} and calculate the $_5F_4$ series on the
right-hand side by \eqref{eq:Dougall},  there is
\begin{align}
\sum_{k=0}^{\infty}\bigg(\frac{-1}{4}\bigg)^k\frac{(\frac{1}{2})_k^3(d)_k(1-d)_k}{(1)_{k}^3(\frac{3}{2}-d)_{k}^2}(1-d+7k-6dk+10k^2)
=\frac{4}{\pi}\frac{\Gamma(\frac{3}{2}-d)^2}{\Gamma(1-d)^2}.
\label{eq:wei-dd}
\end{align}

Realize that the series in \eqref{eq:wei-dd} is uniformly convergent
for $d\in \mathbb{C}$. Apply $\mathcal{D}_{d}$ on both sides of
\eqref{eq:wei-dd} to derive
\begin{align*}
&\sum_{k=0}^{\infty}\bigg(\frac{-1}{4}\bigg)^k\frac{(\frac{1}{2})_k^3(d)_k(1-d)_k}{(1)_{k}^3(\frac{3}{2}-d)_{k}^2}(1-d+7k-6dk+10k^2)
\notag\\[1mm]
&\quad\times\Big\{H_{k}(d-1)-H_{k}(-d)+2H_{k}(\tfrac{1}{2}-d)\Big\}
\notag\\[1mm]
&\:+\sum_{k=0}^{\infty}\bigg(\frac{-1}{4}\bigg)^k\frac{(\frac{1}{2})_k^3(d)_k(1-d)_k}{(1)_{k}^3(\frac{3}{2}-d)_{k}^2}(-1-6k)
\notag\\[1mm]
&\:\:=\frac{8}{\pi}\frac{\Gamma(\frac{3}{2}-d)^2}{\Gamma(1-d)^2}
\Big\{\psi(1-d)-\psi(\tfrac{3}{2}-d)\Big\}.
\end{align*}
Fixing $d=\frac{1}{2}$ in the last equation and utilizing
\eqref{digamma-a} and \eqref{digamma-b}, we discover
\eqref{eq:wei-b} after some simplification.
\end{proof}

Thirdly, we plan to prove Theorem \ref{thm-c}.

\begin{proof}[{\bf{Proof of Theorem \ref{thm-c}}}]
Take $(a,c,d,e)=(\frac{1}{2},\frac{1}{2},\frac{1}{2},\frac{1}{2})$
in \eqref{equation-bb} and compute the $_5F_4$ series on the
right-hand side by \eqref{eq:Dougall} to detect
\begin{align}
\sum_{k=0}^{\infty}(-1)^k\frac{(\frac{1}{2})_k^3(1-b)_k^3}{(1)_{k}^3(\frac{3}{2}-b)_{k}(\frac{3}{2}-b)_{2k}}B_k(b)
=\frac{2}{\pi}\frac{\Gamma(\frac{3}{2}-b)^2}{\Gamma(1-b)^2},
\label{eq:wei-ff}
\end{align}
where
\begin{align*}
B_k(b) &=2k(1-b+2k)+\frac{2(1+2k)(1-b+k)^2}{3-2b+4k}.
\end{align*}

It is not difficult to understand that the series in
\eqref{eq:wei-ff} is uniformly convergent for $b\in \mathbb{C}$.
Apply $\mathcal{D}_{b}$ on both sides of \eqref{eq:wei-ff} to gain
\begin{align*}
&\sum_{k=0}^{\infty}(-1)^k\frac{(\frac{1}{2})_k^3(1-b)_k^3}{(1)_{k}^3(\frac{3}{2}-b)_{k}(\frac{3}{2}-b)_{2k}}B_k(b)
\notag\\[1mm]
&\quad\times\Big\{H_{k}(\tfrac{1}{2}-b)+H_{2k}(\tfrac{1}{2}-b)-3H_{k}(-b)\Big\}
\notag\\[1mm]
&\:+\sum_{k=0}^{\infty}(-1)^k\frac{(\frac{1}{2})_k^3(1-b)_k^3}{(1)_{k}^3(\frac{3}{2}-b)_{k}(\frac{3}{2}-b)_{2k}}\mathcal{D}_{b}B_k(b)
\notag\\[1mm]
&\:\:=\frac{4}{\pi}\frac{\Gamma(\frac{3}{2}-b)^2}{\Gamma(1-b)^2}
\Big\{\psi(1-b)-\psi(\tfrac{3}{2}-b)\Big\}.
\end{align*}
Letting $b=\frac{1}{2}$ in the last equation and employing
\eqref{digamma-a} and \eqref{digamma-b}, there holds
\begin{align}
\sum_{k=0}^{\infty}\frac{\binom{2k}{k}^5}{(-2^{12})^k}\Big\{(20k^2+8k+1)(5H_{2k}-4H_k)+14k+3\Big\}=\frac{32\log2}{\pi^2}.
 \label{eq:wei-gg}
\end{align}
Therefore, the linear combination of \eqref{eq:wei-b} with
\eqref{eq:wei-gg} produces \eqref{eq:wei-c}.
\end{proof}
\section{Proof of Theorems
\ref{thm-d}-\ref{thm-g}}

Above all, we shall prove Theorem \ref{thm-d}.

\begin{proof}[{\bf{Proof of Theorem \ref{thm-d}}}]
Recall a transformation formula for hypergeometric series (cf.
\cite[Theorem 14]{Chu-b}):
\begin{align}
&\sum_{k=0}^{\infty}\frac{(c)_k(e)_k(1+a-b-c)_k(1+a-b-e)_{k}(1+a-c-d)_{k}(1+a-d-e)_{k}}{(1+a-c)_{k}(1+a-e)_{k}}
\notag\\[1mm]
&\quad\times\frac{(1+a-b-d)_{2k}\,\omega_k(a,b,c,d,e)}{(1+a-b)_{2k}(1+a-d)_{2k}(1+2a-b-c-d-e)_{2k}}
\notag\\[1mm]
&=\sum_{k=0}^{\infty}(a+2k)\frac{(b)_k(c)_k(d)_k(e)_k}{(1+a-b)_{k}(1+a-c)_{k}(1+a-d)_{k}(1+a-e)_{k}},
\label{equation-aaa}
\end{align}
where  $\mathfrak{R}(1+2a-b-c-d-e)>0$ and
\begin{align*}
\omega_k(a,b,c,d,e)&=\frac{(1+2a-b-c-d+3k)(a-e+k)}{1+2a-b-c-d-e+2k}
\\[1mm]
&\quad+\frac{(e+k)(1+a-b-c+k)}{(1+a-b+2k)(1+a-d+2k)}
\\[1mm]
&\quad\times\frac{(1+a-c-d+k)(1+a-b-d+2k)(2+2a-b-d-e+3k)}{(1+2a-b-c-d-e+2k)(2+2a-b-c-d-e+2k)}.
\end{align*}
Select $(a,b,d,e)=(\frac{1}{2},\frac{1}{2},\frac{1}{2},\frac{1}{2})$
in \eqref{equation-aaa} and evaluate the $_5F_4$ series on the
right-hand side by \eqref{eq:Dougall} to  get
\begin{align}
\sum_{k=0}^{\infty}\bigg(\frac{1}{16}\bigg)^k\frac{(\frac{1}{2})_k(\frac{1}{2})_{2k}(c)_k(1-c)_k^2}{(1)_{k}^3(\frac{3}{2}-c)_{k}(\frac{3}{2}-c)_{2k}}
E_k(c) =\frac{2}{\pi}\frac{\Gamma(\frac{3}{2}-c)^2}{\Gamma(1-c)^2},
\label{eq:wei-ccc}
\end{align}
where
\begin{align*}
E_k(c) &=2k(1-c+3k)+\frac{3(1+4k)(1-c+k)^2}{2(3-2c+4k)}.
\end{align*}

Notice that the series in \eqref{eq:wei-ccc} is uniformly convergent
for $c\in \mathbb{C}$. Apply $\mathcal{D}_{c}$ on both sides of
\eqref{eq:wei-ccc} to obtain

\begin{align*}
&\sum_{k=0}^{\infty}\bigg(\frac{1}{16}\bigg)^k\frac{(\frac{1}{2})_k(\frac{1}{2})_{2k}(c)_k(1-c)_k^2}{(1)_{k}^3(\frac{3}{2}-c)_{k}(\frac{3}{2}-c)_{2k}}
E_k(c)
\notag\\[1mm]
&\quad\times\Big\{H_{k}(c-1)+H_{k}(\tfrac{1}{2}-c)-2H_{k}(-c)-H_{2k}(\tfrac{1}{2}-c)\Big\}
\notag\\[1mm]
&\:+\sum_{k=0}^{\infty}\bigg(\frac{1}{16}\bigg)^k\frac{(\frac{1}{2})_k(\frac{1}{2})_{2k}(c)_k(1-c)_k^2}{(1)_{k}^3(\frac{3}{2}-c)_{k}(\frac{3}{2}-c)_{2k}}\mathcal{D}_{c}E_k(c)
\notag\\[1mm]
&\:\:=\frac{4}{\pi}\frac{\Gamma(\frac{3}{2}-c)^2}{\Gamma(1-c)^2}
\Big\{\psi(1-c)-\psi(\tfrac{3}{2}-c)\Big\}.
\end{align*}
Setting $c=\frac{1}{2}$ in the last equation and using
\eqref{digamma-a} and \eqref{digamma-b}, we catch hold of
\eqref{eq:wei-d}.
\end{proof}

Then we are ready to prove Theorem \ref{thm-e}.

\begin{proof}[{\bf{Proof of Theorem \ref{thm-e}}}]
Choose $(a,b,c,e)=(\frac{1}{2},\frac{1}{2},\frac{1}{2},\frac{1}{2})$
in \eqref{equation-aaa} and calculate the $_5F_4$ series on the
right-hand side by \eqref{eq:Dougall} to deduce
\begin{align}
\sum_{k=0}^{\infty}\bigg(\frac{1}{4}\bigg)^k\frac{(\frac{1}{2})_k^3(1-d)_k^2(1-d)_{2k}}{(1)_{k}^3(\frac{3}{2}-d)_{2k}^2}
F_k(d) =\frac{1}{\pi}\frac{\Gamma(\frac{3}{2}-d)^2}{\Gamma(1-d)^2},
\label{eq:wei-ddd}
\end{align}
where
\begin{align*}
F_k(d)=k(1-d+3k)+\frac{(1+2k)(1-d+k)(1-d+2k)(2-d+3k)}{(3-2d+4k)^2}.
\end{align*}

Realize that the series in \eqref{eq:wei-ddd} is uniformly
convergent for $d\in \mathbb{C}$. Apply $\mathcal{D}_{d}$ on both
sides of \eqref{eq:wei-ddd} to gain
\begin{align*}
&\sum_{k=0}^{\infty}\bigg(\frac{1}{4}\bigg)^k\frac{(\frac{1}{2})_k^3(1-d)_k^2(1-d)_{2k}}{(1)_{k}^3(\frac{3}{2}-d)_{2k}}
F_k(d)
\notag\\[1mm]
&\quad\times\Big\{2H_{2k}(\tfrac{1}{2}-d)-H_{2k}(-d)-2H_{k}(-d)\Big\}
\notag\\[1mm]
&\:+\sum_{k=0}^{\infty}\bigg(\frac{1}{4}\bigg)^k\frac{(\frac{1}{2})_k^3(1-d)_k^2(1-d)_{2k}}{(1)_{k}^3(\frac{3}{2}-d)_{2k}}
\mathcal{D}_{d}F_k(d)
\notag\\[1mm]
&\:\:=\frac{2}{\pi}\frac{\Gamma(\frac{3}{2}-d)^2}{\Gamma(1-d)^2}
\Big\{\psi(1-d)-\psi(\tfrac{3}{2}-d)\Big\}.
\end{align*}
Fixing $d=\frac{1}{2}$ in the last equation and utilizing
\eqref{digamma-a} and \eqref{digamma-b}, we have
\begin{align}
\sum_{k=0}^{\infty}\frac{\binom{2k}{k}^4\binom{4k}{2k}}{2^{16k}}\Big\{(120k^2+34k+3)(2H_{4k}+H_{2k}-2H_k)+52k+8\Big\}=\frac{128\log2}{\pi^2}.
 \label{eq:wei-eee}
\end{align}
So the difference of \eqref{eq:wei-d} and \eqref{eq:wei-eee}
engenders \eqref{eq:wei-e}.
\end{proof}

Subsequently, we want to prove Theorem \ref{thm-f}.

\begin{proof}[{\bf{Proof of Theorem \ref{thm-f}}}]
Recall another transformation formula for hypergeometric series (cf.
\cite[Theorem 27]{Chu-b}):
\begin{align}
&\sum_{k=0}^{\infty}\frac{(b)_k(c)_k(d)_k(1+a-b-c)_k(1+a-b-d)_{k}(1+a-c-d)_{k}}{(b+e-a)_{k}(c+e-a)_{k}(d+e-a)_{k}}
\notag\\[1mm]
&\quad\times\frac{(e)_{3k}}{(1+a-b)_{2k}(1+a-c)_{2k}(1+a-d)_{2k}}\mu_k(a,b,c,d,e)
\notag\\[1mm]
&+\frac{\Gamma(1+a-b)\Gamma(1+a-c)\Gamma(1+a-d)\Gamma(1+a-e)}{\Gamma(b)\Gamma(c)\Gamma(d)\Gamma(e)}
\notag\\[1mm]
&\quad\times\frac{\Gamma(b+e-a)\Gamma(c+e-a)\Gamma(d+e-a)\Gamma(1+2a-b-c-d-e)}{\Gamma(1+a-b-c)\Gamma(1+a-b-d)\Gamma(1+a-c-d)}
\notag\\[1mm]
&=\sum_{k=0}^{\infty}(a+2k)\frac{(b)_k(c)_k(d)_k(e)_k}{(1+a-b)_{k}(1+a-c)_{k}(1+a-d)_{k}(1+a-e)_{k}},
\label{equation-fff}
\end{align}
where
\begin{align*}
\mu_k(a,b,c,d,e)&=\frac{(a-c+2k)(a-e)}{a-c-e-k}
-\frac{(c+k)(e+3k)(a-e)(1+a-b-d+k)}{(1+a-d+2k)(a-b-e-k)(a-c-e-k)}
\\[1mm]
&\quad+\frac{(b+k)(c+k)(e+3k)(1+e+3k)}{(1+a-b+2k)(1+a-c+2k)(1+a-d+2k)}
\\[1mm]
&\qquad\times\frac{(a-e)(1+a-b-d+k)(1+a-c-d+k)}{(a-b-e-k)(a-c-e-k)(a-d-e-k)}.
\end{align*}
We may manipulate the
$(a,b,c,e)=(\frac{1}{2},\frac{1}{2},\frac{1}{2},0)$ case of
\eqref{equation-fff} as
\begin{align}
\sum_{k=0}^{\infty}\frac{(\frac{1}{2})_k^3(1)_{3k}(d)_{k}(1-d)_{k}^2}{(1)_{k}^2(1)_{2k}^2(\frac{1}{2}+d)_{k}(\frac{3}{2}-d)_{2k}}
U_k(d) =\frac{1-2d}{\pi\cot(\pi d)}, \label{eq:wei-ggg}
\end{align}
where
\begin{align*}
U_k(d)=\frac{3d-3d^2+(6+10d-11d^2)k+(37-10d)k^2+37k^3}{3(3-2d+4k)}.
\end{align*}

It is clear that the series in \eqref{eq:wei-ggg} is uniformly
convergent for $d\in \mathbb{C}$. Apply $\mathcal{D}_{d}$ on both
sides of \eqref{eq:wei-ggg} to find
\begin{align*}
&\sum_{k=0}^{\infty}\frac{(\frac{1}{2})_k^3(1)_{3k}(d)_{k}(1-d)_{k}^2}{(1)_{k}^2(1)_{2k}^2(\frac{1}{2}+d)_{k}(\frac{3}{2}-d)_{2k}}
U_k(d)
\notag\\[1mm]
&\quad\times\Big\{H_{k}(d-1)-2H_{k}(-d)+H_{2k}(\tfrac{1}{2}-d)-H_{k}(d-\tfrac{1}{2})\Big\}
\end{align*}

\begin{align*}
&\:+\sum_{k=0}^{\infty}\frac{(\frac{1}{2})_k^3(1)_{3k}(d)_{k}(1-d)_{k}^2}{(1)_{k}^2(1)_{2k}^2(\frac{1}{2}+d)_{k}(\frac{3}{2}-d)_{2k}}
 \mathcal{D}_{d}U_k(d)
\notag\\[1mm]
&\:\:=\frac{\pi(1-2d)\csc(\pi d)^2-2\cot(\pi d)}{\pi\cot(\pi d)^2}.
\end{align*}
Taking $d=\frac{1}{2}$ in the last equation, we are led to
\eqref{eq:wei-f} after some simplification.
\end{proof}

Finally, we shall certify Theorem \ref{thm-g}.

\begin{proof}
Let $(a,b,c,d)=(\frac{1}{2},\frac{1}{2},\frac{1}{2},\frac{1}{2})$ in
\eqref{equation-fff} and compute the $_5F_4$ series on the
right-hand side by \eqref{eq:Dougall} to detect
\begin{align}
\sum_{k=0}^{\infty}\frac{(\frac{1}{2})_k^6(1+e)_{3k}}{(1)_{2k}^3(1+e)_{k}^3}
V_k(e) =\frac{\Gamma(1+e)^2\Gamma(\frac{1}{2}-e)^2}{\pi^3}-
\frac{e^2}{\pi}\frac{\Gamma(\frac{1}{2}-e)^2}{\Gamma(1-e)^2},
\label{eq:wei-hhh}
\end{align}
where
\begin{align*}
V_k(e)=\frac{e+5e^2+(3+24e+42e^2)k+(27+108e)k^2+74k^3}{16(e+3k)}.
\end{align*}

It is evident that the series in \eqref{eq:wei-ggg} is uniformly
convergent for $e\in \mathbb{C}$. Apply $\mathcal{D}_{e}$ on both
sides of \eqref{eq:wei-hhh} to discover
\begin{align*}
&\sum_{k=0}^{\infty}\frac{(\frac{1}{2})_k^6(1+e)_{3k}}{(1)_{2k}^3(1+e)_{k}^3}
V_k(e)\Big\{H_{3k}(e)-3H_{k}(e)\Big\}
\notag\\[1mm]
&\:\:+\sum_{k=0}^{\infty}\frac{(\frac{1}{2})_k^6(1+e)_{3k}}{(1)_{2k}^3(1+e)_{k}^3}
 \mathcal{D}_{e}V_k(e)
\notag\\[1mm]
&\:\:=\frac{2\Gamma(1+e)^2\Gamma(\frac{1}{2}-e)^2}{\pi^3}\Big\{\psi(1+e)-\psi(\tfrac{1}{2}-e)\Big\}
\notag\\[1mm]
&\quad\:\:-\frac{2e^2}{\pi}\frac{\Gamma(\frac{1}{2}-e)^2}{\Gamma(1-e)^2}\Big\{\psi(1-e)-\psi(\tfrac{1}{2}-e)\Big\}.
\end{align*}
 Thanks to \eqref{digamma-a} and \eqref{digamma-b}, the $e=0$ case of the last equation
 can be stated as
\begin{align}
\sum_{k=0}^{\infty}\frac{\binom{2k}{k}^4\binom{3k}{k}}{2^{12k}}\Big\{3(74k^2+27k+3)(H_{3k}-3H_k)+250k+45\Big\}=\frac{576\log2}{\pi^2}.
 \label{eq:wei-iii}
\end{align}
Thus the linear combination of \eqref{eq:wei-f} with
\eqref{eq:wei-iii} brings out \eqref{eq:wei-g}.
\end{proof}


\end{document}